\documentclass[12pt, notitlepage]{amsart}
\usepackage{latexsym, amsfonts, amsmath, amssymb,bbm, amsthm}

\newtheorem{lemma}{Lemma}
\newtheorem{proposition}{Proposition}
\newtheorem{theorem}{Theorem}
\newtheorem{corollary}{Corollary}

\newcommand{\assign}{:=}

\newcommand{\mathd}{\mathrm{d}}
\newcommand{\mathi}{\mathrm{i}}

\newcommand{\um}{-}

\usepackage{graphicx}
\usepackage{mathrsfs}

\begin{document}

\title{The $4.36$-th moment of the Riemann zeta-function}

\begin{abstract}
  Conditionally on the Riemann Hypothesis we obtain bounds of the correct
  order of magnitude for the $2 k$-th moment of the Riemann zeta-function for
  all positive real $k < 2.181$. 
  This provides for the first time an upper
  bound of the correct order of magnitude for some $k > 2$; the case of $k =
  2$ corresponds to a classical result of Ingham {\cite{Ingh}}. We prove our
  result by establishing a connection between moments with $k > 2$ and the
  so-called ``twisted fourth moment''. This allows us to appeal to a recent
  result of Hughes and Young {\cite{Hughes}}. Furthermore we obtain a
  point-wise bound for $| \zeta (\tfrac{1}{2} + \mathrm{i} t) |^{2 r}$ (with $0 < r <
  1$) that can be regarded as a multiplicative analogue of Selberg's bound for
  $S (T)$ {\cite{Selberg}} . We also establish asymptotic formulae for moments
  ($k < 2.181$) slightly off the half-line.
\end{abstract}

\author{Maksym Radziwi\l\l}
\address{Department of Mathematics \\ Stanford University\\
450 Serra Mall, Bldg. 380\\
Stanford, CA 94305-2125}
\email{maksym@stanford.edu}
\thanks{The author is partially supported by a NSERC PGS-D award}

\maketitle

\section{Introduction.}

\noindent An important problem in analytic number theory is to gain an understanding of
the moments of the Riemann zeta-function,
\begin{equation}
  \label{Moment} M_k (T) := \int_T^{2 T} \left| \zeta \left(\tfrac{1}{2} +
  \mathrm{i} t \right) \right|^{2 k} \mathrm{d} t . 
\end{equation}
It is conjectured that $M_k (T) \sim C_k \cdot T (\log T)^{k^2}$ for all $k >
0$. Recent evidence for this conjecture comes from many sources. Classical
analytical arguments for $k = 6$ and $k = 8$ (see {\cite{Ghosh}} and
{\cite{Gonek}}), models based on random matrix theory (see
{\cite{KeatingSnaith}} and {\cite{KeatingGonek}}), and considerations of
multiple Dirichlet series (see {\cite{Goldfeld}}) ... have all led to the same
conjecture for $M_k (T)$.

Following the work of Hardy-Littlewood {\cite{HarLil}}, and Ingham
{\cite{Ingh}}, this asymptotic formulae is known for $k = 0, 1, 2$. Lower
bounds of the correct order of magnitude have been established for $2 k \in
\mathbbm{N}$ by Ramachandra {\cite{Ram}}, for $2 k \in \mathbbm{Q}$ by
Heath-Brown {\cite{Brown}}, and recently for all $2 k > 0$ by the author and
Soundararajan {\cite{SoundRadz1}}.

Less is known when it comes to upper bounds for $M_k (T)$. Assuming the
Riemann Hypothesis Heath-Brown {\cite{Brown}} established an upper bound of
the correct order of magnitude for $M_k (T)$ for $0 \leqslant k < 2$. For $k >
2$, we have Soundararajan's {\cite{SoundII}} upper bound,
\begin{equation}
  M_k (T) \ll_{k, \varepsilon} T (\log T)^{k^2 + \varepsilon}
\end{equation}
(Here $\varepsilon = O (1 /\ensuremath{\operatorname{logloglog}}T)$ as
remarked by Ivi\'c {\cite{Ivic}}). While these results come close to the
expected order of growth of $M_k (T)$, they miss the correct order of
magnitude of $M_k (T)$ for all $k > 2$. In this paper, we develop a new method
which obtains on the Riemann Hypothesis an upper bound of the right order of
magnitude for the $2 k$-th moment of $\zeta (s)$ in the range $2 < k < 2 + 
\tfrac{2}{11}$.

\begin{theorem}
  \label{Thm1}Assume the Riemann Hypothesis. If $0 < k < 2 + \tfrac{2}{11}$, then,
  \[ M_k (T) \ll_k T (\log T)^{k^2}. \]
\end{theorem}

Our method also provides a new way to obtain Heath-Brown's {\cite{Brown}}
bounds in the range $0 < k < 2$. One surprising feature (which is also present
in {\cite{Brown}}) is that the implicit constant $c (k)$ in our bounds for
$M_k (T)$ tends to infinity as $k \rightarrow 2^-$ even though we have
Ingham's asymptotic formulae for $k = 2$ and now also good upper bounds for $2
< k < 2 + \tfrac{2}{11}$. In particular $c (k)$ stays bounded as $k \rightarrow 2^+$
from the left.

As a quick application of Theorem \ref{Thm1}, 
we obtain asymptotic formulae for
moments $(k < 2 + \tfrac{2}{11})$ slightly off the half-line. 
These are known for all
$\sigma \geqslant \tfrac{1}{2}$ for $k = 1, 2$ by the classical work of
Hardy-Littlewood {\cite{HarLil}} and Ingham {\cite{Ingh}} and following
Soundararajan {\cite{SoundII}} we can also establish asymptotic formulae for
all $k > 0$ when $\sigma \geqslant \tfrac{1}{2} + c_k \cdot
\ensuremath{\operatorname{loglog}}T / \log T$. However, as shown below, we can
do better in terms of $\sigma$ for moments with $k < 2 + \tfrac{2}{11}$. \

\begin{corollary}
  \label{Thm2}Assume the Riemann Hypothesis. 
  Let $\psi (T)$ be such
  that $\psi (T) \rightarrow \infty$ arbitrarily slowly and $\psi (T) = o
  (\log T)$. Set $\sigma = \tfrac 12 + \psi / \log T$. 
  If $0 < k < 2 + \tfrac{2}{11}$, then,
  \[ \int_T^{2 T} \left| \zeta \left( \sigma +
  \mathrm{i} t \right) \right|^{2 k} \mathrm{d} t =
  \sum_{n \geq 1} \frac{d_k(n)^2}{n^{2\sigma}} 
  \cdot \left( 1 + O \left( e^{- 
     \psi / 100} \right) \right). \]
\end{corollary}

Returning to Theorem 1, a natural line of attack is to consider the so-called
``twisted fourth moment''
\begin{equation}
  \label{Twisted} \int_T^{2 T} | \zeta (\tfrac{1}{2} + \mathrm{i} t) |^4 \cdot |A (\tfrac{1}{2} + \mathrm{i} t) |^2 \mathrm{d} t,
\end{equation}
with $A (\cdot)$ a Dirichlet polynomial. The expression in (\ref{Twisted}) was
first considered by Iwaniec in {\cite{Iwaniec}} and then Deshouillers and
Iwaniec in {\cite{IwaniecDesI}} and {\cite{DesIwan2}}. They obtained a bound
$T^{1 + \varepsilon}$ for Dirichlet polynomials $A (s)$ of length 
$T^{1/5 -
\varepsilon}$ and with bounded coefficients (say). Their work has been
substantially refined by Watt {\cite{Watt}}, who allows for Dirichlet
polynomials of length up to $T^{1/4 - \varepsilon}$. An asymptotic formulae
for (\ref{Twisted}) has been recently developed by Hughes and Young
{\cite{Hughes}} (see also related work by Jara {\cite{Jara}}) and this forms
one of the key ingredients in our proof.

Until now the known results about (\ref{Twisted}) had no consequence on
(\ref{Moment}) even on the assumption of the Riemann Hypothesis; the main
difficulty being the absence of a bound for $| \zeta (\tfrac{1}{2} + \mathrm{i} t)
|^{2 r}$ ($0 < r < 1$) in terms of a short Dirichlet polynomial.
Our proof of Theorem 1 provides such a connection for 
all $2 < k < 3$. 
\begin{corollary} \label{cor2}
Assume the Riemann Hypothesis. Let $0 < r < 1$. Then,
\begin{equation}
  \label{Connection} M_{2 + r}(T) \ll \int_{T}^{2T} | \zeta(\tfrac{1}{2}
+ \mathi t)|^{4} \cdot \left | \sum_{n \leqslant x} \frac{d_r(n)W(n)}
{n^{\tfrac{1}{2} + \mathi t}} \right |^{2} \mathd t
\end{equation}
where $x = T^{r / 2 + 2 \delta} < T$ and $\delta > 0$ is an arbitrary but
fixed positive real number, and $W (\cdot)$ is a continuous smoothing, defined as
\begin{equation}
  \label{W} W (n) := \begin{cases}
  1 & \text{ if } n \leqslant y \\
  \frac{\log (x/n)}{\log (x/y)} & \text{ if } y < n \leqslant x \\
  0 & \text{ if } n > x
  \end{cases}
  \text{ where } \left\{\begin{array}{l}
    x = T^{r / 2 + 2 \delta}\\
    y = T^{r / 2 + \delta}
  \end{array}\right.
\end{equation}
The implicit constant in (\ref{Connection}) depends at most on $r$ and
$\delta$.
\end{corollary}
The main idea in our proof is to obtain bounds for $| \zeta (s) |^{2
r}$ $(0 < r < 1)$ in terms of a Dirichlet polynomial and we do this in a way
remniscent of Selberg's work on $S (T)$ {\cite{Selberg}}. More precisely, we
write
\begin{equation}
  \zeta (s)^r = \sum_{n \leqslant x} \frac{d_r (n) W (n)}{n^s} + O \left(
  \frac{y^{\tfrac{1}{2} - \sigma}}{\log (x / y)} \int_{- \infty}^{\infty} \frac{|
  \zeta (\tfrac{1}{2} + \mathrm{i} t + \mathrm{i} v) |^r}{(\sigma - \tfrac{1}{2})^2 + v^2}
  \mathd v \right )
\end{equation}
and upon averaging over $T \leqslant t \leqslant 2 T$ we bound the
contribution from the error term in terms of the moment itself times a small
constant $< 1$. We provide below a variant of this idea which might of
independent interest.

\begin{proposition}
  \label{Prop3}Assume the Riemann Hypothesis. Let $0 < r \leqslant 1$ and $0 <
  \delta < 1$ such that $r / 2 + 3 \delta < 1$ be given. Then, for $T \leqslant t \leqslant 2 T$ with $T
  \geqslant T_0$ (and $T_0$ a large absolute constant),
  \begin{equation}
    \label{ineq} | \zeta (\tfrac{1}{2} + \mathrm{i} t) |^{2 r} \leqslant 3
    e^{12 r / \delta} \int_{\mathbbm{R}} \left| \sum_{n \leqslant x} \frac{d_r
    (n) W (n)}{n^{2 \sigma_0 - \tfrac{1}{2} + \mathrm{i} t + \mathrm{i} v}} \right|^2
    \cdot \frac{(1/ \pi)(\sigma_0 - \tfrac{1}{2}) \mathrm{d} v}{(\sigma_0 - \tfrac{1}{2})^2 +
    v^2} + O \left( \frac{1}{T} \right)
  \end{equation}
  where $x = T^{r / 2 + 2 \delta}$ and $\sigma_0 = \tfrac{1}{2} + \frac{(4 /
  \delta)}{\log T}$. 
\end{proposition}

Integrating the above inequality over $T \leqslant t \leqslant 2 T$ we recover
Heath-Brown's {\cite{Brown}} upper bound $M_r (T) \ll T (\log T)^{r^2}$ for
all $0 \leqslant r < 2$. Furthermore, Proposition \ref{Prop3} allows us to 
generalize (\ref{Connection}) to a product of two (or more!) 
Dirichlet polynomials and also to obtain bounds of the correct order of 
magnitude for the ``twisted $r$-th moment'' $\int_{T}^{2T} | \zeta(\tfrac 12
+ \mathi t)|^{2r} \cdot |A(\tfrac 12 + \mathi t)|^{2} \mathd t$
for $r < 2$ and Dirichlet polynomials of length $< T^{1 - r/2 - \varepsilon}$ and 
for $2 < r < 2 + \tfrac{2}{11}$ and Dirichlet polynomials of
length $T^{1/11 - r/2 - \varepsilon}$ when combined with Hughes and Young's
work \cite{Hughes}.
 
The appearance of $2 \sigma_0 - \tfrac{1}{2} + \mathrm{i} t$
in (\ref{ineq}) is
inessential : we can further bound the right-hand side by
$$
  \leqslant C (\varepsilon) \int_{- \infty}^{\infty} \left| \sum_{n
  \leqslant x} \frac{d_r (n) W (n)}{n^{\tfrac{1}{2} + \mathrm{i} t + \mathrm{i} v}}
  \right|^2 \cdot \frac{(1 / \pi) (\sigma_0 - \tfrac{1}{2}) \mathrm{d} v}{(\sigma_0 -
  \tfrac{1}{2})^2 + v^2},
$$
if we wish to do so. It is also worth noticing that the error term 
$O (1 / T)$ in (\ref{ineq}) is always smaller than the main term. 

Proposition \ref{Prop3} can be regarded as a multiplicative analogue of Soundararajan's
{\cite{SoundII}} upper bound for $\log | \zeta (\tfrac{1}{2} + \mathrm{i} t) |$,
\[ \log | \zeta (\tfrac{1}{2} + \mathrm{i} t) | \leqslant \mathfrak{R} \sum_{n
   \leqslant x} \frac{\Lambda (n) / \log n}{n^{\tfrac{1}{2} + \tfrac{\lambda}{\log x}
   + \mathrm{i} t}} \cdot \frac{\log (x / n)}{\log x} + \frac{1 + \lambda}{2}
   \cdot \frac{\log T}{\log x} + O \left( \frac{1}{\log x} \right) \]
where $T \leqslant t \leqslant 2 T$. Note however, the greater degree of
flexibility in the choice of parameters $1 \leqslant x \leqslant T^2$ and
$\lambda \geqslant \lambda_0 = 0.4912 \ldots$ in the inequality above.

Besides Selberg's work on $S (T)$ our proof is also inspired by ideas from
Ramachandra's and Heath-Brown's work. Our proof adapts to the case of central
moments of $L$-functions (giving for example an alternative proof of the
result in {\cite{BrownDir}} in the range $0 < k < 2$) and we plan to return to
this subject on a later occasion.

{\textbf{Acknowledgment. }}I would like to thank my supervisor Prof.
Soundararajan for his advice and a simplification of my original proof, and
Prof. Matt Young for a number of insightful remarks.

{\textbf{Notation. }}We denote by $\varepsilon$ an arbitrarily small but fixed
constant, not necessarily the same from line to line and by $T_0$ a large
absolute constant. We will write $s = \sigma + \mathrm{i} t$.

\section{Key ideas}

The relevance of the following lemma to moments of the Riemann zeta-function
(on the assumption of the Riemann Hypothesis) was first pointed out by Soundararajan
in {\cite{SoundII}}.

\begin{lemma}
  \label{Lemma4}Let $\xi (s) := s (s - 1) \pi^{- s / 2} \Gamma (s / 2) \zeta
  (s)$. The Riemann Hypothesis holds if and only if for every fixed $t$, $| \xi (\sigma +
  \mathrm{i} t) |$ is an increasing function of $\sigma$ for $\sigma \geqslant
  \tfrac{1}{2}$ . In particular on the Riemann
  Hypothesis, uniformly in $\sigma' \geqslant \sigma \geqslant \tfrac{1}{2}$,
  \[ | \zeta (\sigma + \mathrm{i} t) |^{2 r} \leqslant T^{r (\sigma' -
     \sigma)} \cdot | \zeta (\sigma' + \mathrm{i} t) |^{2 r}, \]
  for all $T \leqslant t \leqslant 2 T$ with $T \geqslant T_0$. 
\end{lemma}

\noindent\textbf{Proof.\ }We have the following Hadamard factorization formula,
\[ \xi (s) = \xi (0) \prod_{\rho} \left( 1 - \frac{s}{\rho} \right) \]
where we group together the zeroes $\rho$ and $\bar{\rho}$ to ensure the
convergence of the product. If the Riemann Hypothesis holds, then
$\Re (\rho) = \tfrac{1}{2}$ for all $\rho$'s, and hence if $\Re s
\geqslant \tfrac{1}{2}$ and $c > 0$, then
\[ | \rho - s| \leqslant | \rho - (s + c) | \]
for any $\rho$, because the distance of $s + c$ from any given point on the
line $\mathfrak{R}z = \tfrac{1}{2}$ is always greater than that of $s$. Thus
\[ \left | \frac{\xi (s)}{\xi (s + c)} \right | \leqslant \prod_{\rho} \left| \frac{\rho -
   s}{\rho - (s + c)} \right| \leqslant 1, \]
as claimed. To establish the converse suppose to the contrary that the Riemann
Hypothesis is false, but $| \xi (\sigma + \mathi t) |$ is still an increasing
function of $\sigma$. Let $\rho = \sigma' + \mathi t$ denote a zero of $\zeta
(s)$ off the half-line. Pick an $\tfrac{1}{2} \leqslant \sigma < \sigma'$ for which $|
\xi (\sigma + \mathi t) | > 0$; such a $\sigma$ exists because $\xi$ is
analytic and not identically zero. Since $| \xi (\sigma + \mathi t) |$ is
increasing $0 < | \xi (\sigma + \mathi t) | \leqslant | \xi (\sigma' + \mathi
t) | = 0$ and this is a contradiction.

Finally, if the Riemann Hypothesis holds, then the inequality $| \zeta
(\sigma + \mathi t) | \leqslant T^{(\sigma' - \sigma)} \cdot | \zeta (\sigma'
+ \mathi t) |$ follows from $| \xi (\sigma + \mathi t) | \leqslant | \xi
(\sigma' + \mathi t) |$ upon unfolding $\xi (s)$ into $s (s - 1) \pi^{- s / 2}
\Gamma (s / 2) \zeta (s)$ and using Stirling's
formula.\hspace*{\fill}$\Box$\medskip

\begin{lemma}
  \label{Lemma5}Assume the Riemann Hypothesis and let $W (n)$ be defined as in
  $( \ref{W})$. Then, for any fixed $r > 0$, uniformly in $\tfrac{1}{2} \leqslant \sigma \leqslant 1$, $T
  \leqslant t \leqslant 2 T$ (with $T \geqslant T_0$) and $1 \leqslant y \leqslant x/2 \leqslant T^{1 -
  \varepsilon}$,
  \begin{eqnarray*}
    \zeta (s)^r & = & \sum_{n \geqslant 1} \frac{d_r (n) W (n)}{n^s} +
    \frac{\theta y^{\tfrac{1}{2} - \sigma}}{\log (x / y)} \cdot \frac{1}{\pi} \int_{-
    \infty}^{\infty} \frac{| \zeta (\tfrac{1}{2} + \mathrm{i} t + \mathrm{i} v)
    |^r}{(\sigma - \tfrac{1}{2})^2 + v^2} \cdot \mathrm{d} v + O \left( \frac{1}{T}
    \right),
  \end{eqnarray*}
  with a $\theta$ of absolute value at most one. 
\end{lemma}

\noindent\textbf{Proof.\ }Consider
\[ \sum_{n \geqslant 1} \frac{d_r (n) W (n)}{n^s} = \frac{1}{2 \pi \mathi}
   \int_{(c)} \zeta (s + w)^r \cdot \frac{x^w - y^w}{w^2 \cdot \log (x / y)}
   \mathd w, \]
with $c =  1 + 1/\log T$.
In order to avoid a branch cut emanating from $w = 1 -s$,  
we take out the little piece corresponding to 
$|\mathfrak{I} (w - s)| \leqslant \varepsilon$ from the integral above,
making an error of $O(1/T)$. 
Shifting contours in the remaining two integrals, we obtain $\zeta(s)^{r}$
from a pole at $w = 0$ and an error of $O(1/T)$ from the integrals
over horizontal lines, since $|\zeta(s)| \ll T^{\varepsilon}$ for
$\sigma \geqslant \tfrac{1}{2}$ by Lindel\"of's Hypothesis.  Thus,
\[ \sum_{n \geqslant 1} \frac{d_r (n) W (n)}{n^s} = \zeta (s)^r + \frac{1}{2
   \pi \mathi} \int_{\mathcal{L}} \zeta (s + w)^r \cdot \frac{(x^w - y^w) \mathd w}{w^2 \cdot \log (x /
   y)} + O \left( \frac{1}{T} \right). \]
where $\mathcal{L}$ corresponds to the line $\mathfrak{R} s = \tfrac{1}{2} 
- \sigma$ without the segment with 
$|\mathfrak{I}(w - s)| \leqslant \varepsilon$. 
Since the integral on the right is in absolute value at most
\[ \frac{\theta y^{\tfrac{1}{2} - \sigma}}{\log (x / y)} \cdot \frac{1}{\pi} \int_{-
   \infty}^{\infty} \frac{| \zeta (\tfrac{1}{2} + \mathi t + \mathi v) |^r}{(\sigma -
   \tfrac{1}{2})^2 + v^2} \cdot \mathd v, \]
the lemma follows. \hspace*{\fill}$\Box$\medskip

\begin{lemma}
  \label{Lemma6}Assume the Riemann Hypothesis. Then uniformly in $\tfrac{1}{2}
  < \sigma \leqslant 1$ and $T \leqslant t \leqslant 2 T$ with $T
  \geqslant T_0$,
  \begin{equation*}
    \int_T^{2 T} | \zeta (\tfrac{1}{2} + \mathrm{i} t) |^4 \cdot \left( \int_{-
    \infty}^{\infty} \frac{| \zeta (\tfrac{1}{2} + \mathrm{i} t + \mathrm{i} v)
    |^r}{(\sigma - \tfrac{1}{2})^2 + v^2} \cdot \mathrm{d} v \right)^2 \mathrm{d} t
    \leqslant
    \frac{2 \pi^2}{(\sigma - \tfrac{1}{2})^2} \int_T^{2 T} |
    \zeta (\tfrac{1}{2} + \mathrm{i} t) |^{2 k} \cdot \mathrm{d} t + O \left( T^{1 -
    \varepsilon} \right),
  \end{equation*}
  where $2 k = 4 + 2 r$. 
\end{lemma}

\noindent\textbf{Proof.\ }The left-hand side is by Cauchy's inequality at most,
\begin{eqnarray*}
  &  & \frac{\pi}{\sigma - \tfrac{1}{2}} \int_T^{2 T} | \zeta (\tfrac{1}{2} + \mathi t) |^4
  \cdot \int_{- \infty}^{\infty} \frac{| \zeta (\tfrac{1}{2} + \mathi t + \mathi v)
  |^{2 r}}{(\sigma - \tfrac{1}{2})^2 + v^2} \mathd v \mathd t\\
  & = & \frac{\pi}{\sigma - \tfrac{1}{2}} \int_{- \infty}^{\infty} \int_T^{2 T} |
  \zeta (\tfrac{1}{2} + \mathi t) |^4 \cdot | \zeta (\tfrac{1}{2} + \mathi t + \mathi v)
  |^{2 r} \mathd t \cdot \frac{\mathd v}{(\sigma - \tfrac{1}{2})^2 + v^2}.
\end{eqnarray*}
The inner integrand is at most $| \zeta (\tfrac{1}{2} + \mathi t) |^{4 + 2 r} + |
\zeta (\tfrac{1}{2} + \mathi t + \mathi v) |^{4 + 2 r}$ . The first term contributes
in total $(\pi / (\sigma - \tfrac{1}{2}))^2 \int_T^{2 T} | \zeta (\tfrac{1}{2} + \mathi t)
|^{4 + 2 r} \mathd t$, while the second makes a contribution of
\[ \leqslant \frac{\pi}{\sigma - \tfrac{1}{2}} \int_{- \infty}^{\infty} \int_T^{2 T}
   | \zeta (\tfrac{1}{2} + \mathi t + \mathi v) |^{4 + 2 r} \mathd t \cdot
   \frac{\mathd v}{(\sigma - \tfrac{1}{2})^2 + v^2}. \]
Consider the integral over $v$. By Lindel\"of's Hypothesis the integral over $t$
is bounded by $T (T + |v|)^{\varepsilon}$ . 
Therefore the terms with $|v| \geqslant
T^{1 - \varepsilon}$ contribute $\ll T^{\varepsilon}$ due to the rapid decay
of the kernel $((\sigma - \tfrac{1}{2})^2 + v^2)^{- 1}$. On the other hand the
contribution from the terms $|v| \leqslant T^{1 - \varepsilon}$ is bounded by
\begin{equation*}
  \int_{- T^{1 - \varepsilon}}^{T^{1 - \varepsilon}} \int_{T - T^{1 -
  \varepsilon}}^{2 T + T^{1 - \varepsilon}} | \zeta (\tfrac{1}{2} + \mathi t) |^{2 k}
  \cdot \mathd t \cdot \frac{\mathd v}{(\sigma - \tfrac{1}{2})^2 + v^2}\\
  \leqslant \frac{\pi}{\sigma - \tfrac{1}{2}} \int_T^{2 T} | \zeta (\tfrac{1}{2} +
  \mathi t) |^{2 k} \mathd t + O \left( T^{1 - \varepsilon} \right).
\end{equation*}
Combining these bounds we obtain the claim.\hspace*{\fill}$\Box$\medskip

\section{The work of Hughes and Young.}

The main result (Theorem 1) in Hughes and Young paper \cite{Hughes} holds
unconditionally. However, we'll aid ourselves in the proof of the next lemma 
by assuming the Riemann Hypothesis. 

\begin{lemma}
  Assume the Riemann Hypothesis. 
  \label{Lemma7}Let $W (n)$ denote the coefficients in $( \ref{W})$. If $0
  \leqslant x \leqslant T^{1 / 11 - \varepsilon}$ then,
  \begin{equation}
    \label{statLemma7} \int_T^{2 T} | \zeta (\tfrac{1}{2} + \mathrm{i} t) |^4 \cdot
    \left| \sum_{n \leqslant x} \frac{d_r (n) W (n)}{n^{\tfrac{1}{2} + \mathrm{i} t}}
    \right|^2 \mathrm{d} t \ll T (\log T)^{(2 + r)^2},
  \end{equation}
  and where the implicit constant depends at most on $\varepsilon$ and $r$.  
\end{lemma}

\noindent\textbf{Proof.\ } We will use freely the notation used in Theorem 1 of
Hughes and Young's paper {\cite{Hughes}}. We pick a smooth function $g$ with
$g (t) = 1$ in $[T ; 2 T]$, $g (t) \leqslant 1$ elsewhere, $g^{(j)} (t) \ll
T^{- j (1 - \varepsilon)}$ for all $j \geqslant 0$ and with support contained
in $[T / 2 ; 3 T]$. By Lemma \ref{Lemma4} 
the left-hand side of (\ref{statLemma7})
is
$$
\ll \int_{T}^{2T} |\zeta(\tfrac 12 + \alpha + \mathi t)|^2 \cdot
|\zeta(\tfrac 12 + \beta + \mathi t)|^{2} \cdot \left | \sum_{n \leq x}
\frac{d_r(n)W(n)}{n^{\tfrac 12 + \mathi t}} \right |^{2} \cdot g(t) \mathd t
$$
where $\alpha = 1 / \log T$ and $\beta = 2 / \log T$. 
By the main result of {\cite{Hughes}}
for $(m,n) = 1$, 
\begin{equation*} 
  \int_{\mathbbm{R}} \left( \frac{m}{n} \right)^{\mathi t} \cdot | \zeta
  (\tfrac{1}{2} + \alpha + \mathi t) |^2 \cdot |\zeta(\tfrac 12
  + \beta + \mathi t)|^{2} g (t) \mathd t\\
  =  
  \int_{\mathbbm{R}} g(t) \mathd t \cdot \frac{Z_{\alpha, \beta, \alpha, \beta, m, n}(0)}{\sqrt{m n}}
  + \ldots + \mathcal{E},
\end{equation*}
where $\mathcal{E} \ll T^{3 / 4 + \varepsilon} \cdot (m n)^{7 / 8}$ and inside
$\ldots$ we omitted a sum of five more additional terms resembling the expression
displayed above 
(the only essential difference being that instead of $Z_{\alpha, \beta, \alpha, \beta, m , n}(0)$ we will have $Z_{\alpha', \beta', \gamma', \delta',m,n}(0)$ with
$\alpha', \beta', \gamma', \delta' \in \{\pm \alpha, \pm \beta\}$. Our choice 
of $\alpha, \beta$ ensures that each of these five
$Z_{\alpha', \beta', \gamma', \delta',m,n}(0)$ is at most $C (\log T)^{4}$ times a 
multiplicative function of $(m,n)$). 

This being said, the left-hand side of (\ref{statLemma7}) is
\begin{align}
  \ll & \int_{\mathbbm{R}} | \zeta (\tfrac{1}{2} + \alpha + \mathi t) |^2  
  \cdot |\zeta(\tfrac 12 + \beta + \mathi t)|^{2} \cdot \left |
  \sum_{n \leqslant x} \frac{d_r (n) W (n)}{n^{\tfrac{1}{2} + \mathi t}} \right|^2
  g (t) \mathd t \nonumber\\
  = & \sum_{m, n \leqslant x} \frac{d_r (n) W (n) d_r (m) W
  (m)}{\sqrt{mn}} \int_{\mathbbm{R}} \left( \frac{m}{n} \right)^{\mathi t}
  \cdot | \zeta (\tfrac{1}{2} + \alpha + \mathi t) |^2 \cdot
  | \zeta ( \tfrac 12 + \beta + \mathi t)|^{2} g (t) \mathd t \nonumber\\
  = &  
  \sum_{m,n \leq x}
  \frac{d_r ( m) d_r ( n) W ( m) W (
  n)}{\sqrt{mn}} \cdot \frac{(m,n)}{\sqrt{mn}} \int_{\mathbbm{R}} g (t) \mathd
  t \cdot Z_{\alpha, \beta, \alpha, \beta, m', n'}(0) + \ldots + \mathcal{E}'
  \label{MainCauchy},
\end{align}
where $m' \assign m / (m,n)$, $n' \assign n / (m,n)$ and
$\mathcal{E}' \ll T^{1 - \varepsilon}$. The function $Z_{\alpha, \beta, \alpha, 
\beta, m' , n'}(0)$ factors into
$$
Z_{\alpha, \beta, \alpha, \beta, m', n'}(0) = A_{\alpha, \beta, \alpha, \beta}(0)
\cdot B_{\alpha, \beta, \alpha, \beta, m', n'}(0) \ll (\log T)^{4} \cdot
|B_{\alpha, \beta, \alpha, \beta, m', n'}(0)|
$$
and $B_{\alpha, \beta, \alpha, \beta, m', n'}(0)$ is multiplicative in 
the two variables $(m,n)$. 
We bound the resulting sum over $m,n$ in (\ref{MainCauchy})
(using $0 \leqslant W(n) \leqslant 1$ for $n \leqslant x$) by
\begin{eqnarray}
  & \leqslant &
  \sum_{m, n \leqslant x} \frac{d_r (m) d_r (n)}{[m,n]} \cdot |B_{\alpha, \beta,
  \alpha, \beta, m', n'} (0) | \nonumber\\
  & \leqslant & \prod_{p \leqslant x} \left( 1 +
  \frac{r}{p} \cdot | \sigma_{\alpha, \beta} (p) | + \frac{r}{p} \cdot |
  \sigma_{\alpha, \beta} (p) | + \frac{r^2}{p} + 
  O \left( \frac{1}{p^2}
  \right) \right)
  \nonumber\\
  & \leqslant & \prod_{p \leqslant x} \left( 1 +
  \frac{r}{p^{1 +\mathfrak{R} \alpha}} + \frac{r}{p^{1 +\mathfrak{R} \beta}} +
  \frac{r}{p^{1 +\mathfrak{R} \alpha}} + \frac{r}{p^{1 +\mathfrak{R} \beta}}
  + \frac{r^2}{p} + O \left( \frac{1}{p^2} \right) \right)  \label{prod}.
\end{eqnarray}
We can bound the sum over $m, n \leqslant x$ by an truncated Euler product
because the summand is a multiplicative function in the two variable $(m,n)$
(i.e $f(m_1 n_1, m_2 n_2) = f(m_1, m_2)f(n_1, n_2)$ for $(m_1 m_2, n_1 n_2) = 1$). 
Since
\[ \sum_{p \leqslant x} \frac{1}{p^{1 + \alpha}} \leqslant \log\log x +
   O (1) \text{ for } \alpha \ll 1 / \log x, \]
the expression in $( \ref{prod})$ is $\ll (\log x)^{r^2 + 4 r} \ll (\log T)^{r^2 + 4r}$. Thus the displayed expression in (\ref{MainCauchy}) is $\ll T (\log T)^{4 + r^2 + 4r} = T (\log T)^{(2 + r)^2}$. The computation of the five additional terms in $\ldots$ is similar and each of them contributes $\ll T (\log T)^{(2 + r)^2}$. We
conclude that (\ref{MainCauchy}) is $\ll T (\log T)^{(2 + r)^2}$ as desired.
\hspace*{\fill}$\Box$

\section{Proof of Theorem \ref{Thm1} (and Corollary \ref{cor2})}

\noindent\textbf{Proof.\ } Theorem \ref{Thm1} is known for $k <
2$, by the work of Heath-Brown {\cite{Brown}}, and for $k = 2$ as a
consequence of Ingham's {\cite{Ingh}} asymptotic formulae. We thus write $k =
2 + r$, with $0 < r < 2 / 11$. By Lemma \ref{Lemma4},
\begin{equation*}
  M_k (T) = \int_T^{2 T} | \zeta (\tfrac{1}{2} + \mathi t) |^4 \cdot | \zeta (
  \tfrac{1}{2} + \mathi t) |^{2 r} \mathd t
  \leqslant T^{r (\sigma - \tfrac{1}{2})} \int_T^{2 T} | \zeta (\tfrac{1}{2} + \mathi
  t) |^4 \cdot | \zeta (\sigma + \mathi t) |^{2 r} \mathd t . 
\end{equation*}
By Lemma \ref{Lemma5} applied to $\zeta (\tfrac{1}{2} + \mathi t)^r$ and the
inequality $|a + b|^2 \leqslant 2| a|^2 + 2| b|^2$ the above is at most
\begin{align}
 \leqslant & 2 T^{r (\sigma - \tfrac{1}{2})} \int_T^{2 T} | \zeta (\tfrac{1}{2} + \mathi t)
   |^4 \cdot \left| \sum_{n \leqslant x} \frac{d_r (n) W (n)}{n^{\tfrac{1}{2} +
   \mathi t}} \right|^2 \mathd t + \hspace*{\fill} \hspace*{\fill}
   \nonumber\\
   \label{ineq_thm1}
    + & 2 T^{r (\sigma - \tfrac{1}{2})} \cdot \frac{y^{1 - 2
   \sigma}}{\log^2 (x / y)} \cdot \frac{1}{\pi^2} \int_T^{2 T} | \zeta (\tfrac{1}{2} + \mathi t) |^4 \cdot
   \left( \int_{- \infty}^{\infty} \frac{| \zeta (\tfrac{1}{2} + \mathi t + \mathi v)
   |^r}{(\sigma - \tfrac{1}{2})^2 + v^2} \mathd v \right)^2 \mathd t,
\end{align}
plus a negligible contribution from $O (1 / T)$, which we henceforth omit. 
By Lemma \ref{Lemma7} the first term in (\ref{ineq_thm1}) is
\begin{equation}
  \label{FirstPf1} \leqslant C \cdot T^{r (\sigma - \tfrac{1}{2})} \cdot T (\log
  T)^{k^2},
\end{equation}
provided that $0 \leqslant x \leqslant T^{1 / 11 - \varepsilon'}$ and with a
finite constant $C = C (\varepsilon' , r)$ depending only on $\varepsilon'$ and $r$
. 
By Lemma \ref{Lemma6} the second term in (\ref{ineq_thm1}) is at most
\begin{equation}
  \label{SecondPf1} 2 T^{r (\sigma - \tfrac{1}{2})} \cdot \frac{y^{1 - 2
  \sigma}}{\log^2 (x / y)} \cdot \frac{2 M_k (T)}{(\sigma - \tfrac{1}{2})^2} +
  O (T).
\end{equation}
One among many possible choices for $x, y$ and $\sigma$ is given by first
fixing an $\delta > 0$ such that $r / 2 + 3 \delta < 1 / 11$, and then setting
$x = T^{r / 2 + 2\delta}$, $y = T^{r / 2 + \delta}$, and $\sigma - 1
/ 2 = (2 / \delta) / \log T$. Note that $x \leqslant T^{1 / 11 -
\delta}$ so that Lemma \ref{Lemma7} is applicable and also that
(\ref{FirstPf1}) is $\leqslant C (\delta, r) e^{2 r / \delta} \cdot T (\log
T)^{k^2}$ by our choice of $\sigma$. In the main term in (\ref{SecondPf1}) we
obtain,
\[ 2 e^{2 r / \delta} \cdot \frac{e^{- 4 (r / 2 + \delta) / \delta}}{\log^2
   (T^{\delta})} \cdot \frac{2 M_k (T)}{(2 / \delta)^2 (\log T)^{- 2}} =
   e^{- 4} \cdot M_k (T) < \frac{M_k (T)}{2} . \]
Therefore $M_k (T) \leqslant C (\delta,r) e^{2 r / \delta} T (\log T)^{k^2} +
M_k (T) / 2 + O(T)$ and hence $M_k (T) \ll T (\log T)^{k^2}$
as claimed.
In order to obtain Corollary \ref{cor2} it suffices not to bound the first term in
(\ref{ineq_thm1}) by (\ref{FirstPf1}) and instead leave it as it is. 
\hspace*{\fill}$\Box$ 

\section{Proof of Corollary \ref{Thm2}}

\noindent\textbf{Proof.\ } Fix an $\delta > 1/100$ (with $\delta < 2/100$)
and in Lemma \ref{Lemma5} choose $x = T^{2 \delta}$, $y = T^{\delta}$, $r = k$ 
and $\sigma = \tfrac{1}{2} + \psi / \log T$. Using the inequality
$|a + e^{-\delta \psi} b|^2 = |a|^2 \cdot ( 1 + O(e^{-\delta \psi}))
+ O(e^{-\delta \psi} |b|^2)$ we get,
\begin{equation*}
 | \zeta (\sigma + \mathi t)^k |^2 = \left| \sum_{n \leqslant x} \frac{d_k
   (n) W (n)}{n^{\sigma + \mathi t}} \right|^2 \left( 1 + O \left( e^{-
   \delta \psi} \right) \right)
   + O \left( \frac{e^{-
   \delta \psi}}{\log^2 (T)} \left( \int_{- \infty}^{\infty} \frac{|
   \zeta (\tfrac{1}{2} + \mathi t + \mathi v) |^k}{(\sigma - \tfrac{1}{2})^2 + v^2} \mathd v
   \right)^2 \right)
\end{equation*}
plus a negligible contribution ($\ll T^{-2}$)
from $O(1/T)$.
Hence,
\begin{align}
  \label{MainPf2} \int_T^{2 T} | \zeta (\sigma + \mathi t) |^{2k} \mathd t = & \int_T^{2 T} \left| \sum_{n \leqslant x} \frac{d_k (n) W
  (n)}{n^{\sigma + \mathi t}} \right|^2 \mathd t \cdot \left ( 1 + O \left ( e^{-\delta \psi} \right ) \right )
  + \nonumber\\
  & + O \left( \frac{e^{- \delta \psi}}{\log^2 (T)} \int_T^{2
  T} \left( \int_{-
  \infty}^{\infty} \frac{| \zeta (\tfrac{1}{2} + \mathi t + \mathi v) |^k}{(\sigma -
  \tfrac{1}{2})^2 + v^2} \mathd v \right)^2 \mathd t \right).
\end{align}
By the same argument as in Lemma
\ref{Lemma6} the error term above is $\ll e^{-\delta \psi}
\int_{T}^{2T} |\zeta(\tfrac 12 + \mathi t)|^{2k} \mathd t + T^{1 - \varepsilon}$
which is $\ll e^{-\delta \psi} T ( \log T)^{k^2} + T^{1 - \varepsilon}$ 
by Theorem 1. On the other hand the main term in (\ref{MainPf2})
is equal to
\[ T \sum_{n \leqslant x} \frac{d_k (n)^2 W (n)^2}{n^{2 \sigma}} \cdot ( 1 +
O(e^{-\delta \psi})) = T \sum_{n \geq 1} \frac{d_k(n)^2}{n^{2\sigma}} \cdot
(1 + O(e^{-\delta \psi})), \]
by an standard calculation. 
Since the main term above is $\asymp T (\log T / \psi)^{k^2}$ it dominates 
the error term. On combining these two estimates the claim follows.  
\hspace*{\fill}$\Box$\medskip

\section{Proof of Proposition \ref{Prop3}.}

\begin{lemma}
  \label{Lemma9}Assume the Riemann Hypothesis. Then, uniformly in 
  $T \leqslant t \leqslant 2T$ (with $T \geqslant T_0$)
  and $\sigma \geqslant \tfrac{1}{2} + 1/\log T$,
  \begin{eqnarray*}
    | \zeta (\tfrac{1}{2} + \mathrm{i} t) |^{2 r} & \leqslant & T^{r (2 \sigma - 1)}
    \int_{- \infty}^{\infty} | \zeta (s + \mathrm{i} v) |^{2 r} \cdot \frac{(1
    / \pi) (\sigma - \tfrac{1}{2}) \mathrm{d} v}{(\sigma - \tfrac{1}{2})^2 + v^2} + O \left(
    \frac{1}{T} \right).
  \end{eqnarray*}
\end{lemma}

\noindent\textbf{Proof.\ } By Cauchy's formula and Lindel\"of's Hypothesis 
(used to
bound the integrals over horizontal lines)
\[ \zeta (s + (\sigma - \tfrac{1}{2}))^{2 r} = \int_{- \infty}^{\infty} \zeta (s + \mathi v)^{2 r} \cdot \frac{(1 / \pi) (\sigma - \tfrac{1}{2}) \mathd
   v}{(\sigma - \tfrac{1}{2})^2 + v^2} + O \left( \frac{1}{T} \right), \]
and by Lemma \ref{Lemma4}, $| \zeta (\tfrac{1}{2} + \mathi t + \mathi v) |^{2 r}
\leqslant T^{r (2 \sigma - 1)} \cdot | \zeta (s + (\sigma - \tfrac{1}{2})) |^{2 r}$.
Combining the two, and taking absolute values we obtain the
claim.\hspace*{\fill}$\Box$\medskip

\begin{lemma}
  \label{Lemma10}
  Let $w_a(v) = \frac{(1/\pi) a}{a^2 + v^2}$ for all $a > 0$. 
  Suppose that $f,g : \mathbbm{R} \rightarrow \mathbbm{R}^+$ are
  two functions in $L^2(\mathbbm{R}; w_a(v) \mathd v)$. If,
  \[ f (x) \leqslant C \cdot g (x) + \frac{1}{4} \int_{- \infty}^{\infty} f (x
     + \mathrm{i} u) w_a (u) \mathrm{d} u, \]
  for some constant $C > 0$, then,
  \[ \int_{- \infty}^{\infty} f (x + \mathrm{i} v)^2 \cdot w_a (v) \mathrm{d}
     v \leqslant 3 C^2 \int_{- \infty}^{\infty} g (x + \mathrm{i} v)^2 \cdot
     w_a (v) \mathrm{d} v. \]
\end{lemma}

\noindent\textbf{Proof.\ }Note that
\begin{eqnarray*}
  f (x)^2 & \leqslant & 2 C^2 \cdot g (x)^2 + \frac{1}{8} \cdot \left( \int_{-
  \infty}^{\infty} f (x + \mathi u) w_a (u) \mathd u \right)^2\\
  & \leqslant & 2 C^2 \cdot g (x)^2 + \frac{1}{8} \cdot \int_{-
  \infty}^{\infty} f (x + \mathi u)^2 w_a (u) \mathd u \cdot \int_{-
  \infty}^{\infty} w_a (u) \mathd u,
\end{eqnarray*}
by the inequality $|a + b|^2 \leqslant 2| a|^2 + 2| b|^2$ and Cauchy-Schwarz.
Also $\int_{- \infty}^{\infty} w_a (u) \mathd u = 1$. In the
above, we let $x \mapsto x + \mathi v$ and integrate with respect to $w_a (v)
\mathd v$ over the whole line $- \infty < v < \infty$. This gives,
\[ \int_{- \infty}^{\infty} f (x + \mathi v)^2 w_a (v) \mathd v \leqslant 2
   C^2 \int_{- \infty}^{\infty} g (x + \mathi v)^2 w_a (v) \mathd v + 
   \frac{1}{8} \int_{- \infty}^{\infty} f (x + \mathi v)^2 \cdot (w_a \ast
   w_a) (v) \mathd v, \]
where $w_a \ast w_a (v) = \int_{- \infty}^{\infty} w_a (u) w_a (u \um v)
\mathd u$. Note that $(w_a \ast w_a) (v) = w_{2 a} (v)$ and that $w_{2 a} (v)
\leqslant 2 w_a (v)$. Whence,
\[ \int_{- \infty}^{\infty} f (x + \mathi v)^2 w_a (v) \mathd v \leqslant 2
   C^2 \int_{- \infty}^{\infty} g (x + \mathi v)^2 w_a (v) \mathd v + \frac{1}{4} \int_{- \infty}^{\infty} f
   (x + \mathi v)^2 w_a (v) \mathd v. \]
Moving the last term to the left-hand side we obtain the claim.
\hspace*{\fill}$\Box$\medskip

\noindent\textbf{Proof of Proposition \ref{Prop3}.\ } In the notation of Lemma
\ref{Lemma10}, let $w (v) \assign w_{\sigma - \tfrac{1}{2}} (v)$. By a simple
modification of Lemma \ref{Lemma5},
\begin{equation}
  \label{Lemma5modif} \zeta (s + \sigma - \tfrac{1}{2})^r = \sum_{n \leqslant x}
  \frac{d_r (n) W (n)}{n^{s + \sigma - \tfrac{1}{2}}} + \frac{\theta y^{\tfrac{1}{2} -
  \sigma}}{\log (x / y)} \cdot \frac{1}{\pi} \int_{- \infty}^{\infty} \frac{|
  \zeta (s + \mathi v) |^r \mathd v}{(\sigma - \tfrac{1}{2})^2 + v^2}
\end{equation}
where $x = T^{r / 2 + 2 \delta}$. Set $y = T^{r / 2 + \delta}$ and $\sigma = 
\tfrac{1}{2} + \frac{4 / \delta}{\log T}$ . By Lemma \ref{Lemma4} and (\ref{Lemma5modif}),
\begin{eqnarray*}
  | \zeta (s) |^r & \leqslant & e^{2 r / \delta} \cdot | \zeta (s + \sigma - 1
  / 2) |^r\\
  & \leqslant & e^{2 r / \delta} \cdot \left| \sum_{n \leqslant x} \frac{d_r
  (n) W (n)}{n^{s + \sigma - \tfrac{1}{2}}} \right| + \frac{1}{4} \int_{-
  \infty}^{\infty} | \zeta (s + \mathi v) |^r \cdot w (v) \mathd v
\end{eqnarray*}
Hence, by Lemma \ref{Lemma10} (with $f(t) = |\zeta(\sigma + \mathi t)|^{r}$
and $g(t) = |\sum_{n \leq x} \frac{d_r(n)W(n)}{n^{2\sigma - 1/2 + \mathi t}}|$),
\[ \int_{- \infty}^{\infty} | \zeta (s + \mathi v) |^{2 r} \cdot w (v) \mathd
   v \leqslant 3 e^{4 r / \delta} \int_{- \infty}^{\infty} \left| \sum_{n
   \leqslant x} \frac{d_r (n) W (n)}{n^{s + \sigma - \tfrac{1}{2} + \mathi v}}
   \right|^2 \cdot w (v) \mathd v \]
By Lemma \ref{Lemma9} the left-hand side is $\geqslant e^{- 8 r / \delta}
\cdot | \zeta (\tfrac{1}{2} + \mathi t) |^{2 r} + O (1 / T)$ and the claim follows. \
\hspace*{\fill}$\Box$

\end{document}